\documentstyle{article}[15pt]
\input amssym.def

\catcode`@=12
 \long\def\@makefntext#1{\noindent #1}
\newskip\tabcentering \tabcentering=1000pt plus 1000pt minus 1000pt

\def\MCH#1#2{\setbox0=\hbox{\raise#1\hbox{#2}}\smash{\box0}}

\def\@evenfoot{}\def\@oddfoot{}

\def\sec#1{\vspace{5mm}\leftline{\bf #1}\vspace{3mm}}

\catcode`@=12

\def\bc{\begin{center}}
\def\ec{\end{center}}

\def\hang{\hangindent\parindent}
\def\textindent#1{\indent\llap{\qquad #1\ \ \enspace}\ignorespaces}
\def\ref{\par\hang\textindent}

\def\a1{(a_1, a_2, \cdots, a_n)}

\def\a{\alpha}

\begin{document}
\thispagestyle{empty}
\vspace*{-3.0truecm}
\noindent
\vspace{1 true cm}
 \bc{\large\bf
Delta shocks  in   the relativistic   full Euler equations for a Chaplygin gas$^{ **}$

\footnotetext{$^{*}$Corresponding author. Tel: +86-0591-83852790.\\
\indent \,\,\,\,\,\,\,\,E-mail address:  zqshao@fzu.edu.cn.\\
\indent \,\,\,\,\,$^{**}$Supported by the National Natural Science
Foundation of China (No. 70371025),  the Scientific Research
Foundation of the
 Ministry of Education of China (No. 02JA790014),
  the Natural Science Foundation of Fujian Province of China   (No.
 2015J01014)   and  the Science and Technology Developmental Foundation
of Fuzhou University (No. 2004-XQ-16).  }}\ec

 \vspace*{0.2 true cm}
\bc{\bf  Zhiqiang  Shao$^{a, *}$\\
{\it $^{a}$Department of Mathematics,  Fuzhou University,  Fuzhou 350002, China}
 }\ec

 \vspace*{2.5 true mm}
\setlength{\unitlength}{1cm}
\begin{picture}(20,0.1)
\put(-0.6,0){\line(1,0){14.5}}
\end{picture}

 \vspace*{2.5 true mm}
\noindent{\small {\small\bf Abstract}

 \vspace*{2.5 true mm}The   relativistic  full Euler equations for a Chaplygin gas
   are studied.  The Riemann problem is solved constructively.
There are two kinds of Riemann solutions, in which one  consists of   three  contact discontinuities and  the other involves
 a  delta
 shock wave on which  both state variables the rest mass density and the proper energy density simultaneously contain the Dirac delta functions. It is quite different from the previous ones on which only one state variable contains the Dirac delta function.
 The formation mechanism, generalized Rankine-Hugoniot relation and  entropy condition are clarified for this type of delta shock wave.
Under  the  generalized  Rankine-Hugoniot
 relation and entropy condition,  the existence and uniqueness of delta shock solutions are also established.

 \vspace*{2.5 true mm}
\noindent{\small {\small\bf MSC: } 35L65;  35L67

 \vspace*{2.5 true mm}
\noindent{\small {\small\bf Keywords:} Relativistic  full Euler equations;  Chaplygin gas;  Riemann  problem;  Delta shock wave

 \vspace*{2.5 true mm}
\setlength{\unitlength}{1cm}
\begin{picture}(20,0.1)
\put(-0.6,0){\line(1,0){14.5}}
\end{picture}



\baselineskip 15pt
 \sec{\Large\bf 1.\quad  Introduction }The relativistic Euler equations of the conservation laws of baryon numbers, momentum and
energy reads (see [8, 27])
$$ \left\{\begin{array}{ll}\Big(\frac{n}{\sqrt{1-v^{2}/c^{2}}}\Big)_{t}+\Big(\frac{nv}{\sqrt{1-v^{2}/c^{2}}}\Big)_x=0,
\\[4mm]\Big(\frac{(p/c^{2}+\rho ) v}{1-v^2/c^2}\Big)_t+\Big(\frac{(p/c^{2}+\rho )v^2}{1-v^2/c^2}+p\Big)_x=0,\\[4mm]
\Big(\frac{(p/c^{2}+\rho )v^2/c^2}{1-v^2/c^2}+\rho\Big)_t+\Big(\frac{(p/c^{2}+\rho ) v}{1-v^2/c^2}\Big)_x=0,\end{array}\right .\eqno{(1.1)}
$$
where $n$, $\rho$,    $p$ and $v$  represent  the rest mass density,   the proper energy density, the pressure and the particle
speed,  respectively, and the constant $c$ is the speed of light.

In his fundamental work of 1948, Taub [42] derived system (1.1)  and then obtained the Hugoniot curve  of the relativistic shocks, and also showed that
$\gamma$, the ratio of specific heats, must be
less than $\frac{5}{3}$.  He gave
a more systematic description of relativistic hydrodynamics in his later work [43]. In 1986, Thompson [44] established several relations on the relativistic shock
curves.  He observed that ``the relativistic shock  equations are much more complicated and do not lend themselves to expressions  that are both simple and general".  Since the  high complexity of the system itself, up to now, there are few results for this system in the literature.
Chen [8] solved the Riemann problem to system (1.1) for the polytropic gas with the equations of state $p=(\gamma-1)c^{2}(\rho-n)$   and $p=kSn^{\gamma}$.
Recently, its vanishing pressure limit problem was studied by Yin and Sheng [48].  Besides,   Geng and Li [15] studied  the non-relativistic global limits of the entropy solutions to the Cauchy problem of  system  (1.1) for the isothermal flow $p=k^{2}\rho$. Ding [14]  proved the
global stability of the strong rarefaction wave to
1-D piston problem of system (1.1) for the polytropic gas with the equations of state $p=(\gamma-1)c^{2}(\rho-n)$   and $p=kSn^{\gamma}$.
Here, we concern with
the  equation of state   is$$p=-\frac{1}{\rho},  \,\,\,\,\,\,\,\,\,\eqno{(1.2)}
$$
which was introduced  by  Chaplygin [5], Tsien  [45] and
von Karman  [19] as  a suitable mathematical approximation for
calculating the lifting force on a wing of an airplane in
aerodynamics.   A gas is called
 a Chaplygin gas  it satisfies  the equation of state (1.2).
The   Chaplygin gas owns a negative pressure
and occurs in certain theories of cosmology. Such a gas has been advertised as a possible model for dark energy [1, 16].

In recent years, astrophysicists have growing interests in the Chaplygin gas dynamics, which
replaces the polytropic equation of state $p(\rho)=\rho^\gamma$  $(\gamma >1)$ (e.g., see [6-7, 18, 37, 40]) with $p(\rho)= -\rho^{-1}$. The typical feature of the Chaplygin
gas dynamics is that the $\delta$-shocks appear in non-zero pressure cases. The Riemann problem was
solved for the nonrelativistic case by Brenier [2] and Serre [35], relativistic case by Cheng and Yang [10],
followed by its vanishing pressure limit problem by  Yin and  Song [49].

In this paper, we are interested in the Riemann problem  for
 (1.1) and (1.2) with
Riemann initial data
$$(n, \rho, v)(0,x)=\left\{\begin{array}{ll}(n_{-},\rho_{-},v_{-}),\,\,\,\,\,\,\,\,\,x<0, \\(n_{+},\rho_{+}, v_{+}),\,\,\,\,\,\,\,\,\,x>0, \end{array}\right .           \eqno{(1.3)}$$
where $\rho_{\pm} >0$,  $n_{\pm} >0$ and $v_{\pm} $ are given constant states. The Riemann problem is a special initial value problem where the initial data consisting of two piecewise constant states are separated by a jump discontinuity at the origin for the one-dimensional hyperbolic systems of conservation laws. It is well known that the Riemann problem is the most fundamental problem in the field of nonlinear hyperbolic conservation laws.  Theories of hyperbolic systems of conservation laws can be found in [3-4, 11, 25, 27, 36, 39] etc.

For the Chaplygin gas,  the considered  relativistic  full Euler equations  possess three linearly degenerate characteristic fields, thus the classical elementary waves   only involves contact discontinuities.  The rarefaction wave curves and the shock wave curves are actually coincided to the so-called contact discontinuities in the state space.  Although the system is much more complicated and the
results are much harder to obtain,  with the help of  the contact discontinuity curves, by the analysis on the physicall relevant region and the method of characteristic analysis,  we   construct Riemann solutions only involving contact discontinuities when $\frac{v_{+}+\frac{1}{\rho_{+}}}{1+\frac{v_{+}}{\rho_{+} c^{2}}}>\frac{v_{-}-\frac{1}{\rho_{-}}}{1-\frac{v_{-}}{\rho_{-} c^{2}}}$.
   However, for the case   $\frac{v_{+}+\frac{1}{\rho_{+}}}{1+\frac{v_{+}}{\rho_{+} c^{2}}}\leq\frac{v_{-}-\frac{1}{\rho_{-}}}{1-\frac{v_{-}}{\rho_{-} c^{2}}}$, we find that the Riemann solution can not be constructed by these classical   contact discontinuities and delta shocks should occur.
In this case,  we rigorously  analyze the formation of mechanism for  delta shock wave with  Dirac delta function in  the state variables the rest mass density and the proper energy density. By the definition of the delta shock wave solution to (1.1) and (1.2) in the sense of distributions, we propose the generalized Rankine-Hugoniot relation and entropy condition for this type of delta shock wave.  Thus both existence and uniqueness of  delta shock wave solutions can be obtained by solving the generalized Rankine-Hugoniot relation under
entropy condition.

In this work, it is proven that
  the delta shock wave with  Dirac
delta function in both state variables the rest mass density and the proper energy density develops in solutions of the   relativistic  full Euler equations for the Chaplygin gas. It is quite different from the previous ones on which only one state variable contains the Dirac delta function. To our knowledge, this type of delta shock wave has not been found in the previous studies on  the relativistic  Euler equations. For related researches of delta
shock waves, we refer the readers to [6, 9-10, 12-13, 17, 20-24, 26, 28-35,  37-38, 41, 46-49]  and the references cited therein  for more details.
For the theory of the delta
shock wave with Dirac delta function in multiple state variables, interested readers may refer to [9, 30-33, 46-47] for further details.
Besides, substantially different from the works [9, 30-31, 46-47], where the delta
shock wave with Dirac delta function in multiple state variables has been found only in some non-strictly hyperbolic systems of conservation laws,  we find this type of delta shock wave in a linearly degenerate and strictly   hyperbolic systems of conservation laws.

The rest of this paper is organized as
 follows. In Sections 2, we first clarify the  physically relevant region where we can persent classical solutions and delta shock waves,
 and deduce the classical   contact discontinuity  curves, then   construct Riemann solutions only involving the classical contact discontinuities. In Section 3,  we   analyze the formation of mechanism for  delta shock wave with  Dirac delta function in  both  the rest mass density and the proper energy density.
  We also  propose the generalized Rankine-Hugoniot and entropy condition for this type of delta shock wave  and then prove   the existence and uniqueness of  delta shock wave solutions under the generalized Rankine-Hugoniot relation and
entropy condition.

\baselineskip 15pt
 \sec{\Large\bf 2.\quad   Preliminaries and classical  Riemann solutions }In this section, we present some preliminary knowledge for
 system (1.1) and  construct classical Riemann solutions   of (1.1)-(1.2) with initial data (1.3).
  The physically relevant region for  solutions  is
$$\Lambda=\bigg\{(n, \rho, v)|n>0, \rho>\frac{1}{c}, |v|<c\bigg\},\eqno{(2.1)}$$that is, the sonic speed $\sqrt{p'(\rho)}$  should be strictly less than
the speed of light (see [8]).

 For any smooth solution, system (1.1) with (1.2)
 can be written in matrix form
$$A\left(\begin{array}{cc}n\\\rho\\v
 \end{array}\right)_t+B\left(\begin{array}{cc}n\\\rho\\v
 \end{array}\right)_x=0,\eqno{(2.2)}
$$
where
$$A=\left(
  \begin{array}{ccc}\frac{1}{\sqrt{1-v^{2}/c^{2}}}& 0 & \frac{nv}{c^{2}(1-v^{2}/c^{2})^{3/2}}\\
    0& \frac{\big(\frac{1}{\rho^{2}c^{2}}+1\big)v}{1-v^{2}/c^{2}} & \frac{\big(-\frac{1}{\rho c^{2}}+\rho\big)(1+v^{2}/c^{2})}{(1-v^{2}/c^{2})^{2}} \\0& \frac{1+\frac{v^{2}}{\rho^{2}c^{4}}}{1-v^{2}/c^{2}} & \frac{\frac{2v}{c^{4}}\big(-\frac{1}{\rho}+\rho c^{2}\big)}{(1-v^{2}/c^{2})^{2}}
  \end{array}\right),\eqno{(2.3)}$$
and
$$ B=\left(
  \begin{array}{ccc}\frac{v}{\sqrt{1-v^{2}/c^{2}}}& 0 & \frac{n}{(1-v^{2}/c^{2})^{3/2}}\\
    0& \frac{v^{2}+1/\rho^{2}}{1-v^{2}/c^{2}} & \frac{2\rho v \big(1-\frac{1}{\rho^{2}c^{2}}\big)}{(1-v^{2}/c^{2})^{2}} \\0& \frac{\big(\frac{1}{\rho^{2}c^{2}}
    +1
   \big )v}{1-v^{2}/c^{2}} & \frac{\big(-\frac{1}{\rho c^{2}}+\rho\big) (1+v^{2}/c^{2})}{(1-v^{2}/c^{2})^{2}}
  \end{array}\right).\eqno{(2.4)}$$
It follows from  (2.3)  and (2.4)   that
$$ A^{-1}B=\left(
  \begin{array}{ccc}v& \frac{\frac{nv}{c^{2}\rho^{2}}(1-v^{2}/c^{2})}{\big(\frac{v^{2}}{\rho^{2} c^{4}}-1\big)\big(-\frac{1}{\rho c^{2}}+\rho\big)} & \frac{n}{1-\frac{v^{2}}{\rho^{2} c^{4}}}\\
    0& \frac{v \big(\frac{1}{\rho^{2}c^{2}}-1\big)}{\frac{v^{2}}{\rho^{2} c^{4}}-1} & \frac{\rho \big(\frac{1}{\rho^{2}c^{2}}-1\big)}{\frac{v^{2}}{\rho^{2} c^{4}}-1} \\0& \frac{-\frac{1}{\rho^{2}}
    (1-v^{2}/c^{2})^{2}}{\big(-\frac{1}{\rho c^{2}}+\rho\big)\big(\frac{v^{2}}{\rho^{2} c^{4}}-1\big)} & \frac{v \big(\frac{1}{\rho^{2}c^{2}}-1\big)}{\frac{v^{2}}{\rho^{2} c^{4}}-1}
  \end{array}\right).\eqno{(2.5)}$$

By (2.5), it is not difficult to see  that   system (1.1) with (1.2) has three real and distinct eigenvalues$$\lambda_{1}=\frac{v-\frac{1}{\rho}}{1-\frac{v}{\rho c^{2}}},\,\,\,\,
\lambda_{2}=v,\,\,\,\,\lambda_{3}=\frac{v+\frac{1}{\rho}}{1+\frac{v}{\rho c^{2}}},
\eqno{(2.6)}
$$
with the  corresponding right  eigenvectors
$$\overrightarrow{r}_1=\bigg(\frac{-n}{\big(\rho-\frac{1}{\rho c^{2}}\big)\big(1-v^{2}/c^{2}\big)},\frac{-1}{1-v^{2}/c^{2}},\frac{1/\rho}{\rho-\frac{1}{\rho c^{2}}}\bigg)^{T},\,\,\overrightarrow{r}_2=(1,0,0 )^T,\,\,
\overrightarrow{r}_3=\bigg(\frac{n}{\big(\rho-\frac{1}{\rho c^{2}}\big)\big(1-v^{2}/c^{2}\big)},\frac{1}{1-v^{2}/c^{2}},\frac{1/\rho}{\rho-\frac{1}{\rho c^{2}}}\bigg)^{T},
\eqno{(2.7)}$$
satisfying
$$ \bigtriangledown\lambda_{i}\cdot \overrightarrow{r_i}\equiv 0\, \,(i=1,2,3).$$
Therefore, system (1.1) with (1.2) is  strictly hyperbolic and fully
 linearly degenerate,  and the associated waves are
 contact discontinuities.

Since system (1.1) with (1.2) and the Riemann data (1.3) are invariant under stretching of coordinates:
$(t, x)\rightarrow (\alpha t, \alpha x)~(\alpha$ is a constant),  we seek the self-similar solution $$(n, \rho,v)(t, x)=(n, \rho,v)(\xi),\,\,\,\,\xi=\frac{x}{t}.$$
 Then  Riemann problem (1.1), (1.2) and (1.3) is reduced to the following boundary value problem of  ordinary differential equations:
$$ \left\{\begin{array}{ll}
   -\xi\Big(\frac{n}{\sqrt{1-v^{2}/c^{2}}}\Big)_{\xi}+\Big(\frac{nv}{\sqrt{1-v^{2}/c^{2}}}\Big)_{\xi}=0,\\
   -\xi\Big(\frac{\big(-\frac{1}{\rho c^{2}}+\rho \big) v}{1-v^2/c^2}\Big)_{\xi}+\Big(\frac{\big(-\frac{1}{\rho c^{2}}+\rho \big)v^2}{1-v^2/c^2}-\frac{1}{\rho}\Big)_{\xi}=0,\\
   -\xi\Big(\frac{\big(-\frac{1}{\rho c^{2}}+\rho \big)v^2/c^2}{1-v^2/c^2}+\rho\Big)_{\xi}+\Big(\frac{\big(-\frac{1}{\rho c^{2}}+\rho \big) v}{1-v^2/c^2}\Big)_{\xi}=0,\end{array}\right .\eqno{(2.8)}
$$
with $(n, \rho,v)(\pm\infty)=(n_{\pm}, \rho_{\pm},v_{\pm}).$\\
\indent
For any smooth solution, system (2.8) can be rewritten as

$$\left(
  \begin{array}{ccc}\frac{v-\xi}{\sqrt{1-v^{2}/c^{2}}}&0
 & \frac{nc^{2}-nv\xi
}{c^{2}(1-v^2/c^2)^{3/2}}  \\0&
    \frac{c^{2}(\frac{1}{\rho^{2}}+v^{2})-\xi v(\frac{1}{\rho^{2}}+c^{2})}{c^{2}-v^{2}} &\frac{(-\frac{1}{\rho}+\rho c^{2})(2vc^{2}-\xi c^{2}-\xi v^{2})}{(c^{2}-v^{2})^{2}}  \\0&
    \frac{(\frac{1}{\rho^{2}}+c^{2})c^{2}v-\xi(c^{4}+\frac{v^{2}}{\rho^{2}})}{c^{2}(c^{2}-v^{2})} &\frac{(-\frac{1}{\rho}+\rho c^{2})(c^{2}+v^{2}-2v\xi)}{(c^{2}-v^{2})^{2}}
  \end{array}\right)\left(\begin{array}{cccc}dn\\ d\rho\\dv
 \end{array}\right)=0.\eqno{(2.9)}$$
It provides either the general solution (constant state)
$$(n, \rho, v)={\mathrm Constant }, $$
or   the singular solutions
$$
\left\{
  \begin{array}{ll}
    \xi=\lambda_{1}=\frac{v-\frac{1}{\rho}}{1-\frac{v}{\rho c^{2}}},\, \\
    d\bigg(\frac{v-\frac{1}{\rho}}{1-\frac{v}{\rho c^{2}}}\bigg)=0,\\\frac{dn}{d\rho} =\frac{n \rho c^{2}}{\rho^{2}c^{2}-1},
  \end{array}
\right.
\eqno{(2.10)}  $$
$$
\left\{
  \begin{array}{ll}
    \xi=\lambda_{2}=v, \\
    d\rho=0,\,\,dv=0,\,\,dn\neq0,
  \end{array}
\right.
\eqno{(2.11)}  $$
$$
\left\{
  \begin{array}{ll}
    \xi=\lambda_{3}=\frac{v+\frac{1}{\rho}}{1+\frac{v}{\rho c^{2}}}, \\
    d\bigg(\frac{v+\frac{1}{\rho}}{1+\frac{v}{\rho c^{2}}}\bigg)=0,\\\frac{dn}{d\rho} =\frac{n\rho c^{2}}{\rho^{2}c^{2}-1}.
  \end{array}
\right.
\eqno{(2.12)}  $$

Integrating (2.10) from $(n_{-}, \rho_{-},v_{-})$  to $(n, \rho,v)$   yields  that
$$\xi=\lambda_{1}=\frac{v-\frac{1}{\rho}}{1-\frac{v}{\rho c^{2}}}
=\frac{v_{-}-\frac{1}{\rho_{-}}}{1-\frac{v_{-}}{\rho_{-} c^{2}}}\,\, {\mathrm and}\,\, \frac{n}{n_{-}}=\sqrt{\frac{(\rho c-1)(\rho c+1)}{(\rho_{-} c-1)(\rho_{-} c+1)}}.
  \eqno{(2.13)}  $$
  Similarly,  we have

  $$\xi=\lambda_{2}=v= v_{-}\, \, \,\,
    \rho=\rho_{-}\,\,\, {\mathrm and}\,\,\,n\neq n_{-},
  \eqno{(2.14)}  $$

$$\xi=\lambda_{3}=\frac{v+\frac{1}{\rho}}{1+\frac{v}{\rho c^{2}}}
=\frac{v_{-}+\frac{1}{\rho_{-}}}{1+\frac{v_{-}}{\rho_{-} c^{2}}}\,\, {\mathrm and}\,\, \frac{n}{n_{-}}=\sqrt{\frac{(\rho c-1)(\rho c+1)}{(\rho_{-} c-1)(\rho_{-} c+1)}}.
  \eqno{(2.15)}  $$

\indent
For a bounded discontinuity at $\xi=\sigma,$ the Rankine-Hugoniot relation holds:

$$ \left\{\begin{array}{ll}
   -\sigma\Big[\frac{n}{\sqrt{1-v^{2}/c^{2}}}\Big]+\Big[\frac{nv}{\sqrt{1-v^{2}/c^{2}}}\Big]=0,\\
   -\sigma\Big[\frac{\big(-\frac{1}{\rho c^{2}}+\rho \big) v}{1-v^2/c^2}\Big]+\Big[\frac{\big(-\frac{1}{\rho c^{2}}+\rho \big)v^2}{1-v^2/c^2}-\frac{1}{\rho}\Big]=0,\\
   -\sigma\Big[\frac{\big(-\frac{1}{\rho c^{2}}+\rho \big)v^2/c^2}{1-v^2/c^2}+\rho\Big]+\Big[\frac{\big(-\frac{1}{\rho c^{2}}+\rho \big) v}{1-v^2/c^2}\Big]=0,\end{array}\right .\eqno{(2.16)}
$$ where $[q]=q -q_{-}$ is the jump of $q$  across the discontinuity  and  $\sigma$ is the velocity of the discontinuity.

Eliminating $\sigma$ in the second and third equations of (2.16), we have
$$\bigg(-\frac{1}{\rho}+\frac{1}{\rho_{-}}\bigg)(\rho-\rho_{-})\bigg(1-\frac{v^{2}}{c^{2}}\bigg)\bigg(1-\frac{v_{-}^{2}}{c^{2}}\bigg)=\bigg(-\frac{1}{\rho c^{2}}+\rho \bigg)\bigg(-\frac{1}{\rho_{-} c^{2}}+\rho_{-} \bigg)(v-v_{-})^{2}.\eqno{(2.17)}$$
Then, from (2.17)  it follows that
$$\frac{(v-v_{-})^{2}}{\bigg(1-\frac{v^{2}}{c^{2}}\bigg)\bigg(1-\frac{v_{-}^{2}}{c^{2}}\bigg)}=\frac{(\rho-\rho_{-})^{2}}
{\bigg(\rho^{2}-\frac{1}{ c^{2}} \bigg)\bigg(\rho_{-}^{2}-\frac{1}{ c^{2}} \bigg)},\eqno{(2.18)}$$
and
$$\bigg(\frac{v-v_{-}}{\frac{vv_{-}}{c^{2}}-1}\bigg)^{2}=\frac{\big(-\frac{1}{\rho}+\frac{1}{\rho_{-}}\big)(\rho -\rho_{-})}{\big(-\frac{1}{\rho c^{2}}+\rho _{-} \big)\big(-\frac{1}{\rho_{-} c^{2}}+\rho \big)}.\eqno{(2.19)}$$
By a simple calculation, it is easy to see that (2.19) is equivalent to
$$\frac{v-v_{-}}{\frac{vv_{-}}{c^{2}}-1}=\pm\bigg(\frac{\rho -\rho_{-}}{\rho\rho _{-} -\frac{1}{c^{2}}} \bigg).\eqno{(2.20)}$$
Thus, we have two cases, namely,

Case 1: $\frac{v-v_{-}}{\frac{vv_{-}}{c^{2}}-1}=\frac{\rho -\rho_{-}}{\rho\rho _{-} -\frac{1}{c^{2}}},$  which gives a 1-shock

 $$S_{1}:\frac{v-\frac{1}{\rho}}{1-\frac{v}{\rho c^{2}}}
=\frac{v_{-}-\frac{1}{\rho_{-}}}{1-\frac{v_{-}}{\rho_{-} c^{2}}}, \,\, {\mathrm with}\,\,p>p_{-},\,\rho>\rho_{-},\,v<v_{-}.
  \eqno{(2.21)}  $$
Hence, from (2.21) and (2.18) it follows that
$$\frac{v-v_{-}}{\sqrt{\Big(1-\frac{v^{2}}{c^{2}}\Big)\Big(1-\frac{v_{-}^{2}}{c^{2}}\Big)}}=-\frac{\rho-\rho_{-}}{\sqrt{\Big(\rho^{2}-\frac{1}{ c^{2}} \Big)\Big(\rho_{-}^{2}-\frac{1}{ c^{2}} \Big)}},\eqno{(2.22)}$$
$$v=\frac{v_{-}-a}{1-\frac{v_{-}a}{c^{2}}},\eqno{(2.23)}$$
$$\sqrt{1-\frac{v^{2}}{c^{2}}}=\frac{\sqrt{c^{2}-a^{2}}}{c\Big(1-\frac{v_{-}a}{c^{2}}\Big)}\sqrt{1-\frac{v_{-}^{2}}{c^{2}}}, \eqno{(2.24)}$$
$$\sqrt{c^{2}-a^{2}}=\frac{c\sqrt{\rho\rho_{-}\big(-\frac{1}{\rho c^{2}}+\rho \big)\big(-\frac{1}{\rho_{-} c^{2}}+\rho_{-} \big)}}{\rho\rho_{-}
-\frac{1}{c^{2}}},\eqno{(2.25)}$$
where  $$a=\frac{\rho-\rho_{-}}{\rho\rho_{-}-\frac{1}{c^{2}}}.$$

Eliminating $\sigma$ in the first   and second  equations of (2.16), we have
$$\frac{v-v_{-}}{\sqrt{\Big(1-\frac{v^{2}}{c^{2}}\Big)\Big(1-\frac{v_{-}^{2}}{c^{2}}\Big)}}\Bigg(\frac{n_{-}v\big(-\frac{1}{\rho c^{2}}+\rho \big)}{\sqrt{1-\frac{v^{2}}{c^{2}}}}-\frac{nv_{-}\big(-\frac{1}{\rho_{-} c^{2}}+\rho_{-} \big)}{\sqrt{1-\frac{v_{-}^{2}}{c^{2}}}}\Bigg)=\Big(-\frac{1}{\rho}+\frac{1}{\rho_{-}}\Big)\Bigg(\frac{n}{\sqrt{1-\frac{v^{2}}{c^{2}} }} -\frac{n_{-}}{\sqrt{1-\frac{v_{-}^{2}}{c^{2}}}}\Bigg).\eqno{(2.26)}$$
Substituting (2.22)-(2.25) into (2.26), we get, after a straightforward calculation that$$\Bigg(\frac{n\rho_{-}(\rho\rho_{-}-\frac{1}{c^{2}})}{\sqrt{\rho\rho_{-}\big(-\frac{1}{\rho c^{2}}+\rho \big)\big(-\frac{1}{\rho_{-} c^{2}}+\rho_{-} \big)}}-\frac{n_{-}(\rho\rho_{-}-\frac{1}{c^{2}})}{-\frac{1}{\rho_{-} c^{2}}+\rho_{-} }\Bigg)\Big(v_{-}-\frac{1}{\rho_{-}}\Big)(\rho-\rho_{-})=0.
  \eqno{(2.27)}  $$

When  $\rho\neq\rho_{-}$, the second part of the left side in the above expression will not be zero if $v_{-}\neq\frac{1}{\rho_{-}}$,   which means
$$\frac{n}{n_{-}}=\sqrt{\frac{(\rho c-1)(\rho c+1)}{(\rho_{-} c-1)(\rho_{-} c+1)}},\,\,\,\,\,\,\,\,\,{\mathrm if }\,\,v_{-}\neq\frac{1}{\rho_{-}}.
  \eqno{(2.28)}  $$

Eliminating $\sigma$ in the first   and third equations of (2.16), we have
$$\frac{nv\Big(-\frac{1}{\rho c^{2}}+\rho\Big)}{\sqrt{1-\frac{v^{2}}{c^{2}}}}+\frac{n_{-}v_{-}\Big(-\frac{1}{\rho_{-} c^{2}}+\rho_{-}\Big)}{\sqrt{1-\frac{v_{-}^{2}}{c^{2}}}}-\frac{n_{-}\Big(-\frac{1}{\rho c^{2}}+\rho\Big)\Big(v-\frac{v_{-}v^{2}}{c^{2}}\Big)}{\Big(1-\frac{v^{2}}{c^{2}}\Big)\sqrt{1-\frac{v_{-}^{2}}{c^{2}}}}$$
$$-\frac{n\Big(-\frac{1}{\rho_{-} c^{2}}+\rho_{-}\Big)\Big(v_{-}-\frac{vv_{-}^{2}}{c^{2}}\Big)}{\Big(1-\frac{v_{-}^{2}}{c^{2}}\Big)\sqrt{1-\frac{v^{2}}{c^{2}}}}=(\rho-\rho_{-})
\Bigg(\frac{nv}{\sqrt{1-\frac{v^{2}}{c^{2}}}}
-\frac{n_{-}v_{-}}{\sqrt{1-\frac{v_{-}^{2}}{c^{2}}}}\Bigg).$$
We rearrange terms to get
$$\frac{nv\Big(-\frac{1}{\rho c^{2}}+\rho_{-}\Big)}{\sqrt{1-\frac{v^{2}}{c^{2}}}}+\frac{n_{-}v_{-}\Big(-\frac{1}{\rho_{-} c^{2}}+\rho\Big)}{\sqrt{1-\frac{v_{-}^{2}}{c^{2}}}}=\frac{1-\frac{vv_{-}}{c^{2}}}{\sqrt{\Big(1-\frac{v^{2}}{c^{2}}\Big)\Big(1-\frac{v_{-}^{2}}{c^{2}}\Big)}}
\Bigg(\frac{n_{-}v\Big(-\frac{1}{\rho c^{2}}+\rho\Big)}{\sqrt{1-\frac{v^{2}}{c^{2}}}}+\frac{nv_{-}\Big(-\frac{1}{\rho_{-} c^{2}}+\rho_{-}\Big)}{\sqrt{1-\frac{v_{-}^{2}}{c^{2}}}}\Bigg).\eqno{(2.29)}  $$
Substituting (2.22)-(2.25)
into (2.29), and  noting that $\frac{v-v_{-}}{\frac{vv_{-}}{c^{2}}-1}=\frac{\rho -\rho_{-}}{\rho\rho _{-} -\frac{1}{c^{2}}}$, we get, after a straightforward calculation that$$\Bigg(\frac{n(\rho\rho_{-}-\frac{1}{c^{2}})}{\sqrt{\rho\rho_{-}\big(-\frac{1}{\rho c^{2}}+\rho \big)\big(-\frac{1}{\rho_{-} c^{2}}+\rho_{-} \big)}}-\frac{n_{-}(\rho\rho_{-}-\frac{1}{c^{2}})}{\rho_{-}\big(-\frac{1}{\rho_{-} c^{2}}+\rho_{-}\big) }\Bigg)\Big(\frac{v_{-}}{\rho_{-}c^{2}}-1\Big)(\rho-\rho_{-})=0.
  \eqno{(2.30)}  $$

When  $\rho\neq\rho_{-}$, the second part of the left side in the above expression will not be zero if $v_{-}\neq\rho_{-}c^{2}$,   which means
$$\frac{n}{n_{-}}=\sqrt{\frac{(\rho c-1)(\rho c+1)}{(\rho_{-} c-1)(\rho_{-} c+1)}},\,\,\,\,{\mathrm if}\,\, v_{-}\neq\rho_{-}c^{2}.
  \eqno{(2.31)}  $$
Because $v_{-}=\rho_{-}c^{2}$  contradicts with $v_{-}=\frac{1}{\rho_{-}}$, the above expression (2.31) together with (2.28) yields that
$$\frac{n}{n_{-}}=\sqrt{\frac{(\rho c-1)(\rho c+1)}{(\rho_{-} c-1)(\rho_{-} c+1)}}.
  \eqno{(2.32)}  $$
This defines the Hugoniot curve  of the relativistic shock.

Substituting (2.23)-(2.25) and (2.32) into the first equation of (2.16), we get, after a straightforward calculation that
$$\sigma=\frac{v_{-}\rho_{-}(\rho-\rho_{-})+a(\frac{1}{c^{2}}-\rho\rho_{-})}{\rho_{-}(\rho-\rho_{-})+\frac{v_{-}a}{c^{2}}(\frac{1}{c^{2}}-\rho\rho_{-})}
=\frac{(\rho-\rho_{-})(v_{-}\rho_{-}-1)}{(\rho-\rho_{-})(\rho_{-}-\frac{v_{-}}{c^{2}})}.\eqno{(2.33)}  $$

When  $\rho\neq\rho_{-}$, from (2.33),  it is easy to find that
$$\sigma
=\frac{v_{-}-\frac{1}{\rho_{-}}}{1-\frac{v_{-}}{\rho_{-} c^{2}}}.\eqno{(2.34)}  $$

When  $\rho=\rho_{-}$, the situation is simple. From (2.23) and (2.16),  we can easily obtain that
$$ \sigma=v= v_{-},\, \, \,\,
    \rho=\rho_{-}\,\,\, {\mathrm and}\,\,\,n\neq n_{-},
  \eqno{(2.35)}  $$

  Case 2: $\frac{v-v_{-}}{\frac{vv_{-}}{c^{2}}-1}=-\frac{\rho -\rho_{-}}{\rho\rho _{-} -\frac{1}{c^{2}}},$  which gives a 3-shock

 $$S_{3}:\frac{v+\frac{1}{\rho}}{1+\frac{v}{\rho c^{2}}}
=\frac{v_{-}+\frac{1}{\rho_{-}}}{1+\frac{v_{-}}{\rho_{-} c^{2}}}, \,\, {\mathrm with}\,\,p<p_{-},\,\rho<\rho_{-},\,v<v_{-}.
  \eqno{(2.36)}  $$
Then, from (2.36) and (2.18) it follows that
$$\frac{v-v_{-}}{\sqrt{\Big(1-\frac{v^{2}}{c^{2}}\Big)\Big(1-\frac{v_{-}^{2}}{c^{2}}\Big)}}=\frac{\rho-\rho_{-}}{\sqrt{\Big(\rho^{2}-\frac{1}{ c^{2}} \Big)\Big(\rho_{-}^{2}-\frac{1}{ c^{2}} \Big)}},\eqno{(2.37)}$$
$$v=\frac{v_{-}+a}{1+\frac{v_{-}a}{c^{2}}},\eqno{(2.38)}$$
$$\sqrt{1-\frac{v^{2}}{c^{2}}}=\frac{\sqrt{c^{2}-a^{2}}}{c\Big(1+\frac{v_{-}a}{c^{2}}\Big)}\sqrt{1-\frac{v_{-}^{2}}{c^{2}}}, \eqno{(2.39)}$$
$$\sqrt{c^{2}-a^{2}}=\frac{c\sqrt{\big(\rho^{2}-\frac{1}{ c^{2}}\big)\big(\rho_{-}^{2}-\frac{1}{ c^{2}} \big)}}{\rho\rho_{-}
-\frac{1}{c^{2}}},\eqno{(2.40)}$$
where  $$a=\frac{\rho-\rho_{-}}{\rho\rho_{-}-\frac{1}{c^{2}}}.$$

Substituting (2.37)-(2.40) into (2.26), we get, after a straightforward calculation that
$$\Bigg(\frac{n(\rho_{-}-\frac{1}{\rho c^{2}})}{\rho_{-}\sqrt{\big(\rho^{2}-\frac{1}{ c^{2}} \big)\big(\rho_{-}^{2}-\frac{1}{ c^{2}} \big)}}-\frac{n_{-}(\rho-\frac{1}{\rho_{-}c^{2}})}{\rho\big(\rho_{-}^{2}-\frac{1}{ c^{2}}\big) }\Bigg)(\rho_{-}v_{-}+1)(\rho-\rho_{-})=0.
  \eqno{(2.41)}  $$

When  $\rho\neq\rho_{-}$, the second part of the left side in the above expression will not be zero if $v_{-}\neq -\frac{1}{\rho_{-}}$,   which means
$$\frac{n}{n_{-}}=\sqrt{\frac{(\rho c-1)(\rho c+1)}{(\rho_{-} c-1)(\rho_{-} c+1)}},\,\,\,\,\,\,\,\,\,{\mathrm if }\,\,v_{-}\neq -\frac{1}{\rho_{-}}.
  \eqno{(2.42)}  $$

Substituting (2.37)-(2.40)
into (2.29), and  noting that $\frac{v-v_{-}}{\frac{vv_{-}}{c^{2}}-1}=-\frac{\rho -\rho_{-}}{\rho\rho _{-} -\frac{1}{c^{2}}}$, we get, after a straightforward calculation that$$\Bigg(\frac{n(\rho\rho_{-}-\frac{1}{c^{2}})}{\sqrt{\big(\rho^{2}-\frac{1}{ c^{2}} \big)\big(\rho_{-}^{2}-\frac{1}{ c^{2}} \big)}}-\frac{n_{-}(\rho\rho_{-}-\frac{1}{c^{2}})}{\rho_{-}^{2}-\frac{1}{ c^{2}} }\Bigg)\Big(\frac{v_{-}}{\rho_{-}c^{2}}+1\Big)(\rho-\rho_{-})=0.
  \eqno{(2.43)}  $$

When  $\rho\neq\rho_{-}$, the second part of the left side in the above expression will not be zero if $v_{-}\neq  -\rho_{-}c^{2}$,   which means
$$\frac{n}{n_{-}}=\sqrt{\frac{(\rho c-1)(\rho c+1)}{(\rho_{-} c-1)(\rho_{-} c+1)}},\,\,\,\,{\mathrm if}\,\, v_{-}\neq -\rho_{-}c^{2}.
  \eqno{(2.44)}  $$
Because $v_{-}=-\rho_{-}c^{2}$  contradicts with $v_{-}=-\frac{1}{\rho_{-}}$, the above expression (2.44) together with (2.42) yields that
$$\frac{n}{n_{-}}=\sqrt{\frac{(\rho c-1)(\rho c+1)}{(\rho_{-} c-1)(\rho_{-} c+1)}}.
  \eqno{(2.45)}  $$
This defines the Hugoniot curve  of the relativistic shock.

Substituting (2.38)-(2.40) and (2.45) into the first equation of (2.16), we get, after a straightforward calculation that
$$\sigma=\frac{v_{-}\rho_{-}(\rho-\rho_{-})+a(\rho\rho_{-}-\frac{1}{c^{2}})}{\rho_{-}(\rho-\rho_{-})+\frac{v_{-}a}{c^{2}}(\rho\rho_{-}-\frac{1}{c^{2}})}
=\frac{(\rho-\rho_{-})(v_{-}\rho_{-}+1)}{(\rho-\rho_{-})(\rho_{-}+\frac{v_{-}}{c^{2}})}.\eqno{(2.46)}  $$

When  $\rho\neq\rho_{-}$, from (2.46),  it is easy to find that
$$\sigma
=\frac{v_{-}+\frac{1}{\rho_{-}}}{1+\frac{v_{-}}{\rho_{-} c^{2}}}.\eqno{(2.47)}  $$

When  $\rho=\rho_{-}$, the situation is simple. From (2.38) and (2.16),  we can easily obtain that
$$ \sigma=v= v_{-},\, \, \,\,
    \rho=\rho_{-}\,\,\, {\mathrm and}\,\,\,n\neq n_{-},
  \eqno{(2.48)}  $$

According to above discussions, there are two types of  shock  curves $S_{i}$  $(i=1,3)$,  which are given by

   $$S_{1}:\left\{
             \begin{array}{ll} \sigma=\frac{v-\frac{1}{\rho}}{1-\frac{v}{\rho c^{2}}}
=\frac{v_{-}-\frac{1}{\rho_{-}}}{1-\frac{v_{-}}{\rho_{-} c^{2}}},\\ \frac{n}{n_{-}}=\sqrt{\frac{(\rho c-1)(\rho c+1)}{(\rho_{-} c-1)(\rho_{-} c+1)}},\,
             \end{array}
           \right.
\eqno{(2.49)}  $$   with$$\,\,\,\,p>p_{-},\,\rho>\rho_{-},\,v<v_{-},
  \eqno{}  $$  and $$S_{3}:\left\{
             \begin{array}{ll} \sigma=\frac{v+\frac{1}{\rho}}{1+\frac{v}{\rho c^{2}}}
=\frac{v_{-}+\frac{1}{\rho_{-}}}{1+\frac{v_{-}}{\rho_{-} c^{2}}},\\ \frac{n}{n_{-}}=\sqrt{\frac{(\rho c-1)(\rho c+1)}{(\rho_{-} c-1)(\rho_{-} c+1)}},\,
             \end{array}
           \right.
\eqno{(2.50)}  $$  with$$\,\,p<p_{-},\,\rho<\rho_{-},\,v<v_{-}.
  \eqno{}  $$

From (2.13),  (2.15) and (2.49)-(2.50), we can find that the rarefaction waves and the shock waves are coincident in the state space, which correspond to contact discontinuities of the first and the third families:

 $$J_{1}:\xi=\sigma=\frac{v-\frac{1}{\rho}}{1-\frac{v}{\rho c^{2}}}
=\frac{v_{-}-\frac{1}{\rho_{-}}}{1-\frac{v_{-}}{\rho_{-} c^{2}}}\,\, {\mathrm and}\,\, \frac{n}{n_{-}}=\sqrt{\frac{(\rho c-1)(\rho c+1)}{(\rho_{-} c-1)(\rho_{-} c+1)}},
  \eqno{(2.51)}  $$

$$J_{3}:\xi=\sigma=\frac{v+\frac{1}{\rho}}{1+\frac{v}{\rho c^{2}}}
=\frac{v_{-}+\frac{1}{\rho_{-}}}{1+\frac{v_{-}}{\rho_{-} c^{2}}}\,\, {\mathrm and}\,\, \frac{n}{n_{-}}=\sqrt{\frac{(\rho c-1)(\rho c+1)}{(\rho_{-} c-1)(\rho_{-} c+1)}}.
  \eqno{(2.52)}  $$

From (2.14),  (2.35) and (2.48), we can find that   there is  a  contact discontinuity  of  the
second family:

  $$J_{2}: \xi=\sigma=v= v_{-},\, \, \,\,
    \rho=\rho_{-}\,\,\, {\mathrm and}\,\,\,n\neq n_{-}.
  \eqno{(2.53)}  $$
\vskip 0.1in\vskip 0.1in
 \noindent{\small {\small\bf Remark 1.} In [8],  the author derived  the Hugoniot curve  of the relativistic shocks  by considering such a coordinate system that the shock  speed is zero. Different from that, in this paper, we obtain the result  generally and ulteriorly give  the analytical formula of relativistic shocks ( also see [15, 48]).
\vskip 0.1in

 In the state  space, starting from a given  state  $(n_{-},\rho_{-},v_{-}),$    we draw the contact discontinuity curves  (2.51) and (2.52)
 for $\rho>\frac{1}{c}$, the projections of  which  onto the $(\rho, v)$-plane are denoted by $J_{1}$ and $J_{3}$,  respectively. So,  $J_{1}$ has the asymptotic line
$v=\frac{v_{-}-\frac{1}{\rho_{-}}}{1-\frac{v_{-}}{\rho_{-} c^{2}}}$ and the singularity point $(\rho, v)=(\frac{1}{c}, c),$ and $J_{3}$ has the asymptotic line
$v=\frac{v_{-}+\frac{1}{\rho_{-}}}{1+\frac{v_{-}}{\rho_{-} c^{2}}}$ and the singularity point $(\rho, v)=(\frac{1}{c}, -c).$
Also, starting from the point $\bigg(n_{-},\rho_{-}, \frac{\rho_{-}v_{-}c^{2}-2c^{2}+\frac{v_{-}}{\rho_{-}}}{\rho_{-}c^{2}-2v_{-}+\frac{1}{\rho_{-}}}\bigg),$
in the state  space, we draw  the contact discontinuity curve (2.52),  the projection of  which  onto the $(\rho, v)$-plane is denoted by $S_{\delta}$.
Then,  $S_{\delta}$  has the asymptotic line
$v=\frac{v_{-}-\frac{1}{\rho_{-}}}{1-\frac{v_{-}}{\rho_{-} c^{2}}}$ and the singularity point $(\rho, v)=(\frac{1}{c}, -c).$
The projections of  these curves   onto the $(\rho, v)$-plane  divide   the $(\rho, v)$-plane  into five regions I, II, III, IV  and V, as shown in Fig. 1.

\hspace{65mm}\setlength{\unitlength}{0.8mm}\begin{picture}(80,66)
\put(-50,0){\line(0,4){6}}
\put(51,0){\line(0,4){6}}\put(-36.5,35){$\cdot$}
\put(-52,2){\vector(0,2){50}}\put(-24,0){\line(0,4){52}}\put(20,-6){$\frac{v_{-}+\frac{1}{\rho_{-}}}{1+\frac{v_{-}}{\rho_{-} c^{2}}}$ }\put(24.2,0){\line(0,4){52}}
 \put(-50,0){\vector(2,0){105}} \put(2,-4){$v_{-}$} \put(-55,49){$\rho$}\put(52,3){$\rho=\frac{1}{c}$}
\put(56,-1){$v$}\put(50,-3){$c$}\put(-53,-3){$-c$}
\put(16,44){$J_{3}$}
\put(-18,45){$J_{1}$}\put(-29,-6){$\frac{v_{-}-\frac{1}{\rho_{-}}}{1-\frac{v_{-}}{\rho_{-} c^{2}}}$}
\put(-16,20){IV }
\put(3,19){$(\rho_{-},
v_{-})$}
\put(2.5,7){II }\put(44,6){$$ }\put(2,32){III}
\put(-45,20){V}\put(-33,40.5){$S_{\delta}$}
\put(21,20){I }\qbezier(50,4)(-14,12)(-22,52)\qbezier(-50,4)(17.5,4)(22,52)
\put(0,3.4){\line(-4,0){50}}\put(0,3.4){\line(4,0){52}}\qbezier(-50,5)(-27.5,4)(-25,52)
\end{picture}
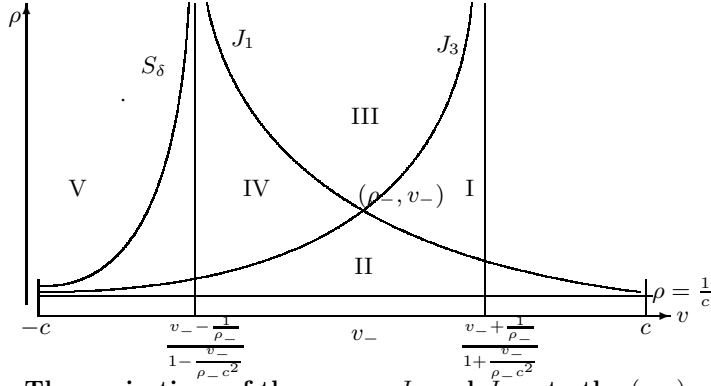
\vspace{0.6mm}  \vskip 0.2in \centerline{\bf Fig. 1.\, The  projections of the curves $J_{1}$ and $J_{3}$  onto the $(\rho, v)$-plane.
   } \vskip 0.1in \indent

\hspace{65mm}\setlength{\unitlength}{0.8mm}\begin{picture}(80,66)\put(0,0){\line(-4,5){25}}
 \put(0,0){\line(1,3){12}}
\put(0,0){\line(3,2){34}}

\put(10,38){$J_{2}$}\put(-30,30){$J_{1}$}
\put(25,8){$(n_{+},\rho_{+},
v_{+})$}
\put(35,22){$J_{3}$}\put(-48,18){$$}\put(-40,8){$(n_{-},\rho_{-},
v_{-})$}

\put(11,26){$(n_{*2},\rho_{*2}, v_{*2})$ }\put(-18,23){$(n_{*1},\rho_{*1}, v_{*1})$ }
\put(-34,34){$$}

 \put(0,0){\vector(0,2){48}}\put(0,0){\vector(0,2){48}}
\put(-45,0){\vector(2,0){98}}
 \put(-3,-4){$O$} \put(-4,45){$t$}
\put(54,-1){$x$}
\end{picture}
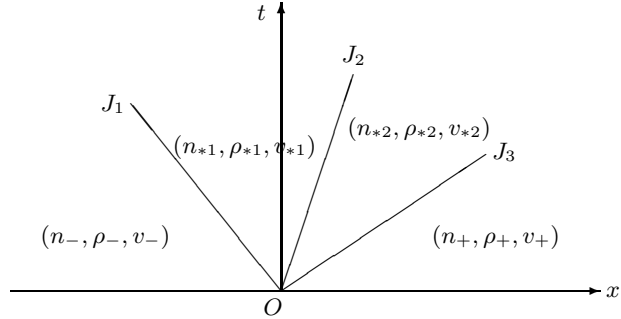
\vspace{0.6mm}  \vskip 0.2in \centerline{\bf Fig. 2.\, The  Riemann solution of the case $\frac{v_{+}+\frac{1}{\rho_{+}}}{1+\frac{v_{+}}{\rho_{+} c^{2}}}>\frac{v_{-}-\frac{1}{\rho_{-}}}{1-\frac{v_{-}}{\rho_{-} c^{2}}}$.} \vskip 0.1in \indent

For any given right state $(n_{+},\rho_{+}, v_{+}),$  according to Fig. 1,  we can construct
Riemann solutions of (1.1), (1.2) and (1.3). When the projection of the state $(n_{+},\rho_{+}, v_{+})$  onto  the $(\rho, v)$-plane lies in  $ I\cup II\cup III\cup IV$, the Riemann problem can be solved in the following way.
On the physically relevant region,  we draw  the contact
discontinuity  curves   $J_{1}(n_{-}, \rho_{-},v_{-})$  and $J_{3}(n_{+}, \rho_{+},v_{+})$.
The projections of these contact
discontinuity  curves  onto  the $(\rho, v)$-plane  have a unique intersection point  $(\rho_{*}, v_{*})$ determined by
$$\left\{\begin{array}{ll}\frac{v_{*}-\frac{1}{\rho_{*}}}{1-\frac{v_{*}}{\rho_{*} c^{2}}}
=\frac{v_{-}-\frac{1}{\rho_{-}}}{1-\frac{v_{-}}{\rho_{-} c^{2}}},\\\frac{v_{*}+\frac{1}{\rho_{*}}}{1+\frac{v_{*}}{\rho_{*} c^{2}}}
=\frac{v_{+}+\frac{1}{\rho_{+}}}{1+\frac{v_{+}}{\rho_{+} c^{2}}}.\end{array}\right .\eqno{(2.54)}$$
Then, we draw the contact
discontinuity  curve   $J_{2}(n_{-}, \rho_{*},v_{*}),$  which intersects the contact
discontinuity  curves   $J_{1}(n_{-}, \rho_{-}, v_{-})$  and $J_{3}(n_{+}, \rho_{+}, v_{+})$ at the unique points $(n_{*1},\rho_{*1}, v_{*1}),$ and $(n_{*2},\rho_{*2}, v_{*2}),$   $\rho_{*1}=\rho_{*2}=\rho_{*}, $  $v_{*1}=v_{*2}=v_{*}. $
So far, we have completely obtained   a solution of  (1.1), (1.2) and (1.3),  as shown in Fig. 2.
Thus,  we  have proved the following result

\vskip 0.1in
\noindent{\small {\small\bf Theorem 2.1.}   For
Riemann problem (1.1),  (1.2)  and (1.3), on the physically relevant region,  under the condition $\frac{v_{+}+\frac{1}{\rho_{+}}}{1+\frac{v_{+}}{\rho_{+} c^{2}}}>\frac{v_{-}-\frac{1}{\rho_{-}}}{1-\frac{v_{-}}{\rho_{-} c^{2}}}$,   there exists   a unique entropy solution,  which  can be expressed  as
$$(n, \rho, v)(t,x)=\left\{
    \begin{array}{ll}
      (n_{-}, \rho_{-}, v_{-}), \,\, \,\,\,\,\,-\infty<x/t<\frac{v_{-}-\frac{1}{\rho_{-}}}{1-\frac{v_{-}}{\rho_{-} c^{2}}},\\(n_{*1},\rho_{*1}, v_{*1}),   \,\,\,\,\, \,\,\,\,\,\frac{v_{-}-\frac{1}{\rho_{-}}}{1-\frac{v_{-}}{\rho_{-} c^{2}}}\leq x/t\leq v_{*1},\\(n_{*2},\rho_{*2}, v_{*2}),   \,\,\,\,\, \,\,\,\,\,v_{*1}<x/t\leq\frac{v_{+}+\frac{1}{\rho_{+}}}{1+\frac{v_{+}}{\rho_{+} c^{2}}},\\
     (n_{+},\rho_{+}, v_{+}),\,\,\,\,\,\,\,\frac{v_{+}+\frac{1}{\rho_{+}}}{1+\frac{v_{+}}{\rho_{+} c^{2}}}<x/t<+\infty,
    \end{array}
  \right.\eqno{(2.55)}$$
 where $$\left\{\begin{array}{ll}\rho_{*1}=\rho_{*2}=\rho_{*}, \\ v_{*1}=v_{*2}=v_{*},\\\frac{v_{*}-\frac{1}{\rho_{*}}}{1-\frac{v_{*}}{\rho_{*} c^{2}}}
=\frac{v_{-}-\frac{1}{\rho_{-}}}{1-\frac{v_{-}}{\rho_{-} c^{2}}},\\\frac{v_{*}+\frac{1}{\rho_{*}}}{1+\frac{v_{*}}{\rho_{*} c^{2}}}
=\frac{v_{+}+\frac{1}{\rho_{+}}}{1+\frac{v_{+}}{\rho_{+} c^{2}}},\,\,\\ n_{*1}=n_{-}\sqrt{\frac{(\rho_{*} c-1)(\rho_{*} c+1)}{(\rho_{-} c-1)(\rho_{-} c+1)}},\\ n_{*2}=n_{+}\sqrt{\frac{(\rho_{*} c-1)(\rho_{*} c+1)}{(\rho_{+} c-1)(\rho_{+} c+1)}}.\end{array}\right .\eqno{(2.56)}$$

\hspace{65mm}\setlength{\unitlength}{0.8mm}\begin{picture}(80,66)
\put(0,0){\line(3,2){40}}\put(0,0){\line(-4,5){25}}\put(0,0){\line(-1,5){8}}
 \put(0,0){\line(1,3){15}}
\put(8,48){$\frac{x}{t}=\frac{v_{+}+\frac{1}{\rho_{+}}}{1+\frac{v_{+}}{\rho_{+} c^{2}}}$}
\put(25,8){$(n_{+},\rho_{+},
v_{+})$}
\put(41,29){$\frac{x}{t}=v_{-}$}
\put(35,40){$\frac{x}{t}= \frac{v_{-}-\frac{1}{\rho_{-}}}{1-\frac{v_{-}}{\rho_{-} c^{2}}}$}\qbezier(3.1,10)(8,13)(9,9)

\put(35,20){$\frac{x}{t}= \frac{v_{-}+\frac{1}{\rho_{-}}}{1+\frac{v_{-}}{\rho_{-} c^{2}}}$}

\put(-48,18){$$}\put(-36,4){$(n_{-},\rho_{-},
v_{-})$}

\put(8.5,16){$\Omega$ }\put(-14,40){$\frac{x}{t}=v_{+}$}

\put(-34,34){$\frac{x}{t}=\frac{v_{+}-\frac{1}{\rho_{+}}}{1-\frac{v_{+}}{\rho_{+} c^{2}}}$}

 \put(0,0){\vector(0,2){48}}\put(0,0){\line(2,1){40}}\put(0,0){\line(1,1){38}}
\put(-45,0){\vector(2,0){98}} \put(-3,-4){$O$} \put(-4,45){$t$}
\put(54,-1){$x$}
\end{picture}
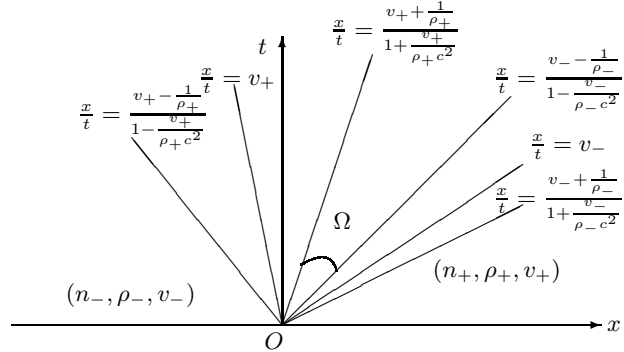
\vspace{0.6mm}  \vskip 0.2in \centerline{\bf Fig. 3.\, The characteristic  lines  from initial data   for   the case
 $\frac{v_{-}-\frac{1}{\rho_{-}}}{1-\frac{v_{-}}{\rho_{-} c^{2}}}\geq \frac{v_{+}+\frac{1}{\rho_{+}}}{1+\frac{v_{+}}{\rho_{+} c^{2}}}$.} \vskip 0.1in \indent

\baselineskip 15pt
 \sec{\Large\bf 3.\quad   Delta shock   solutions }
 In this section, we construct the Riemann solutions   of (1.1)-(1.2) with initial data (1.3)  when  the projection of the state $(n_{+},\rho_{+}, v_{+})$  onto  the $(\rho, v)$-plane lies in  $V$, namely,
 $$\frac{v_{-}-\frac{1}{\rho_{-}}}{1-\frac{v_{-}}{\rho_{-} c^{2}}}\geq \frac{v_{+}+\frac{1}{\rho_{+}}}{1+\frac{v_{+}}{\rho_{+} c^{2}}}.\eqno{(3.1)}$$
 At this moment, the linearly degenerate characteristic  lines from initial data will overlap
in a domain $\Omega=\{(t,x)|\frac{v_{+}+\frac{1}{\rho_{+}}}{1+\frac{v_{+}}{\rho_{+} c^{2}}}t\leq x\leq\frac{v_{-}-\frac{1}{\rho_{-}}}{1-\frac{v_{-}}{\rho_{-} c^{2}}}t, 0\leq t<+\infty\}$ shown in Fig. 3.
So, the singularity
  must
develop  in $\Omega$.  It is easy to know that the singularity is
impossible to be a jump with finite amplitudes because the
Rankine-Hugoniot relation is not satisfied on the bounded jump.

To analyze the singularity, we first study the special case  $\frac{v_{-}-\frac{1}{\rho_{-}}}{1-\frac{v_{-}}{\rho_{-} c^{2}}}= \frac{v_{+}+\frac{1}{\rho_{+}}}{1+\frac{v_{+}}{\rho_{+} c^{2}}},$  let us consider the limit of the solution $ (n,\rho,v)(\xi)$ when
 $n_{-}, $  $\rho_{-}, $  $v_{-}, $  $n_{+}$  and $\rho_{+} $  are fixed, $ \frac{v_{+}+\frac{1}{\rho_{+}}}{1+\frac{v_{+}}{\rho_{+} c^{2}}}\rightarrow \frac{v_{-}-\frac{1}{\rho_{-}}}{1-\frac{v_{-}}{\rho_{-} c^{2}}}+0$.
 When $ \frac{v_{+}+\frac{1}{\rho_{+}}}{1+\frac{v_{+}}{\rho_{+} c^{2}}}> \frac{v_{-}-\frac{1}{\rho_{-}}}{1-\frac{v_{-}}{\rho_{-} c^{2}}}$, the solution
 is given by (2.55),   where $\rho_{*}$ and $ v_{*}$ satisfy
$$\left\{\begin{array}{ll}\frac{c^{2}(v_{*}-\frac{1}{\rho_{*}})}{c^{2}-\frac{v_{*}}{\rho_{*} }}
=\frac{c^{2}(v_{-}-\frac{1}{\rho_{-}})}{c^{2}-\frac{v_{-}}{\rho_{-} }},\\\frac{c^{2}(v_{*}+\frac{1}{\rho_{*}})}{c^{2}+\frac{v_{*}}{\rho_{*} }}
=\frac{c^{2}(v_{+}+\frac{1}{\rho_{+}})}{c^{2}+\frac{v_{+}}{\rho_{+} }}.\end{array}\right .\eqno{(3.2)}$$
We can employ (3.2) and calculate to obtain
$$\rho_{*}=\frac{c^{2}-ab+\sqrt{c^{4}+a^{2}b^{2}-c^{2}(a^{2}+b^{2})}}{c^{2}(b-a)},\eqno{(3.3)}$$
where
$$a=\frac{v_{-}-\frac{1}{\rho_{-}}}{1-\frac{v_{-}}{\rho_{-} c^{2}}} \,\,  { \mathrm \,\, and}\,\,\,\,b = \frac{v_{+}+\frac{1}{\rho_{+}}}{1+\frac{v_{+}}{\rho_{+} c^{2}}}.$$
Therefore, as $ \frac{v_{+}+\frac{1}{\rho_{+}}}{1+\frac{v_{+}}{\rho_{+} c^{2}}}\rightarrow \frac{v_{-}-\frac{1}{\rho_{-}}}{1-\frac{v_{-}}{\rho_{-} c^{2}}}+0$,
namely, $b \rightarrow a^{+}$, the combination of (3.2)-(3.3) and (2.56) yields
$$\left\{\begin{array}{ll}\rho_{*1}=\rho_{*2}=\rho_{*}\rightarrow +\infty, n_{*1}\rightarrow +\infty,n_{*2}\rightarrow +\infty,\\ v_{*1}=v_{*2}=v_{*}\rightarrow
\frac{v_{-}-\frac{1}{\rho_{-}}}{1-\frac{v_{-}}{\rho_{-} c^{2}}}.\end{array}\right .\eqno{(3.4)}$$
These show that   contact
discontinuities  $J_{1}$, $J_{2}$  and $J_{3}$ coincide to form a new  type of nonlinear hyperbolic wave.
Now let us calculate the total quantities of $n$,  $\rho$  and $v$   between $J_{1}$  and   $J_{3}$  as $n_{-}, $  $\rho_{-}, $  $v_{-}, $  $n_{+}$  and $\rho_{+} $  are fixed, $ \frac{v_{+}+\frac{1}{\rho_{+}}}{1+\frac{v_{+}}{\rho_{+} c^{2}}}\rightarrow \frac{v_{-}-\frac{1}{\rho_{-}}}{1-\frac{v_{-}}{\rho_{-} c^{2}}}+0$,
$$\lim\limits_{b\rightarrow a^{+}}\int_{a}^{b}\rho(\xi)d\xi=\lim\limits_{b\rightarrow a^{+}}\int_{a}^{b}\rho_{*}d\xi=\lim\limits_{b\rightarrow a^{+}}\frac{c^{2}-ab+\sqrt{c^{4}+a^{2}b^{2}-c^{2}(a^{2}+b^{2})}}{c^{2}}=\frac{2(c^{2}-a^{2})}{c^{2}}\neq 0,\eqno{(3.5)}$$
$$\lim\limits_{b\rightarrow a^{+}}\int_{a}^{b}v(\xi)d\xi=\lim\limits_{b\rightarrow a^{+}}\int_{a}^{b}v_{*}d\xi=\lim\limits_{b\rightarrow a^{+}}v_{*}(b-a)= 0,\eqno{(3.6)}$$
and
$$\lim\limits_{b\rightarrow a^{+}}\int_{a}^{b}n(\xi)d\xi=\lim\limits_{b\rightarrow a^{+}}\bigg(\int_{a}^{v_{*}}n_{-}\sqrt{\frac{c^{2}\rho_{*}^{2}-1}{c^{2}\rho_{-}^{2}-1}}d\xi
+\int^{b}_{v_{*}}n_{+}\sqrt{\frac{c^{2}\rho_{*}^{2}-1}{c^{2}\rho_{+}^{2}-1}}d\xi\bigg)$$
$$=\lim\limits_{b\rightarrow a^{+}}\bigg(\frac{(c^{2}-v_{*}a)n_{-}}{\rho_{*}c^{2}}\sqrt{\frac{c^{2}\rho_{*}^{2}-1}{c^{2}\rho_{-}^{2}-1}}
+\frac{(c^{2}-v_{*}b)n_{+}}{\rho_{*}c^{2}}\sqrt{\frac{c^{2}\rho_{*}^{2}-1}{c^{2}\rho_{+}^{2}-1}}\bigg)$$
$$=\frac{c^{2}-a^{2}}{c^{2}} \bigg(\frac{n_{-}}{\sqrt{\rho_{-}^{2}-\frac{1}{c^{2}}}}+\frac{n_{+}}{\sqrt{\rho_{+}^{2}-\frac{1}{c^{2}}}}\bigg)\neq0.\eqno{(3.7)}$$
Hence, (3.5) and (3.7) show that $\rho(\xi)$ and $n(\xi)$ have the same singularity as a weighted  Dirac delta function at $\xi=  \frac{v_{-}-\frac{1}{\rho_{-}}}{1-\frac{v_{-}}{\rho_{-} c^{2}}}, $  while (3.6) implies that $v(\xi)$  has a bounded variation.
Thus, such a type of nonlinear hyperbolic wave of (1.1) and (1.2) is
  a delta shock wave with a weighted Dirac
delta function in both $n$ and $\rho $,  denoted by $S_{\delta}$.  It is quite different from the previous ones on which only one state variable contains the Dirac delta function. To our knowledge, this type of delta shock wave has not been found in the previous studies on  the relativistic  Euler equations.
Moreover, for $S_{\delta}$ in this case,  the inequality
$$\lambda_{1}(n_{+},\rho_{+},
v_{+})<\lambda_{2}(n_{+},\rho_{+},
v_{+})<\lambda_{3}(n_{+},\rho_{+},
v_{+})=\sigma=\lambda_{1}(n_{-},\rho_{-},
v_{-})<\lambda_{2}(n_{-},\rho_{-},
v_{-})<\lambda_{3}(n_{-},\rho_{-},
v_{-})$$
holds, where $\sigma=\frac{v_{-}-\frac{1}{\rho_{-}}}{1-\frac{v_{-}}{\rho_{-} c^{2}}}$ is the propagation speed of  $S_{\delta}$.  It means that none of the six
characteristic  lines on both sides of $S_{\delta}$ is outgoing with respect to $S_{\delta}$.

By the above analysis, for   the case
 $\frac{v_{-}-\frac{1}{\rho_{-}}}{1-\frac{v_{-}}{\rho_{-} c^{2}}}\geq \frac{v_{+}+\frac{1}{\rho_{+}}}{1+\frac{v_{+}}{\rho_{+} c^{2}}}$, the delta shock wave  solution
 containing a  Dirac
delta function in both $n$ and $\rho$  will be considered. Thus  we first introduce  three definitions   as  follows.
 \vskip 0.1in

 \vskip 0.1in
\noindent{\small {\small\bf Definition 3.1.} A triple  $(n, \rho, v)$  constitutes a solution of (1.1)  in the sense of distributions if it satisfies
$$ \left\{
     \begin{array}{ll}
       \int_{0}^{+\infty}\int_{-\infty}^{+\infty}\Big(\Big(\frac{n}{\sqrt{1-v^{2}/c^{2}}}\Big)\phi_{t}+\Big(\frac{nv}{\sqrt{1-v^{2}/c^{2}}}\Big)\phi_{x}\Big)dxdt=0,\\
      \int_{0}^{+\infty}\int_{-\infty}^{+\infty}\Big(\Big(\frac{(p/c^{2}+\rho ) v}{1-v^2/c^2}\Big)\phi_{t}+\Big(\frac{(p/c^{2}+\rho )v^2}{1-v^2/c^2}+p\Big)\phi_{x}\Big)dxdt=0,\\\int_{0}^{+\infty}\int_{-\infty}^{+\infty}\Big(\frac{(p/c^{2}+\rho )v^2/c^2}{1-v^2/c^2}+\rho\Big) \phi_{t}+\Big(\frac{(p/c^{2}+\rho ) v}{1-v^2/c^2}\Big)\phi_{x}\Big)dxdt=0,
     \end{array}
   \right.
\eqno{(3.8)}$$
for all test functions $\phi\in C^{\infty}_{0}(R^{+}\times R^{1}).$\\
\vskip 0.1in

\noindent{\small {\small\bf Definition 3.2.}
 A two-dimensional weighted delta function $w(s)\delta_{L}$   supported on a smooth curve $L$ parameterized as
 $t=t(s)$,  $x=x(s)$
 $ (c\leq s\leq d)$ is defined by
$$ \langle w(s)\delta_{L},\phi\rangle=\int_{c}^{d}w(s)\phi(t(s),x(s))ds, \eqno{(3.9)}$$
for all test functions $\phi\in C^{\infty}_{0}(R^{2}).$ \\

\vskip 0.1in

\noindent{\small {\small\bf Definition 3.3.}
 A triple distribution $(n, \rho, v)$ is called a delta shock wave  solution of (1.1)   if
it is represented in the form
$$(n, \rho,v)(t, x)=\left\{
                  \begin{array}{ll}
                    (n_{l}, \rho_{l}, v_{l})(t,x), & \hbox{$x<x(t) $,} \\
                    (h(t)\delta(x-x( t)),w(t)\delta(x-x( t)),   v_{\delta}(t)), & \hbox{$x=x(t)$,} \\
                    (n_{r}, \rho_{r}, v_{r})(t,x), & \hbox{$x>x(t)$}
                  \end{array}
                \right.
\eqno{(3.10)}$$
and satisfies  Definition 3.1, where    $(n_{l}, \rho_{l}, v_{l})(t,x)$  and   $(n_{r}, \rho_{r}, v_{r})(t,x)$   are
 piecewise smooth  bounded solutions of (1.1).

  \vskip 0.1in
\indent

With Definitions 3.1-3.3,  we seek  a delta shock wave solution with the discontinuity $x=x(t)$ of (1.1) with (1.2) in
the form
$$(n, \rho,v)(t, x)=\left\{
                  \begin{array}{ll}
                    (n_{-}, \rho_{-}, v_{-}), & \hbox{$x<x(t) $,} \\
                    (h(t)\delta(x-x( t)),w(t)\delta(x-x( t)),   v_{\delta}(t)), & \hbox{$x=x(t)$,} \\
                    (n_{+}, \rho_{+}, v_{+}), & \hbox{$x>x(t)$,}
                  \end{array}
                \right.
\eqno{(3.11)}$$
where   $x(t), h(t),$  $w(t)\in C^{1}[0, +\infty)$,   $\delta(\cdot)$ is the standard Dirac measure
supported on the curve $x=x(t)$,   and $h(t), w(t)$  are the weights of the delta shock wave on the state variables $n, \rho$, respectively.
Similar to [2, 10, 32], we define $\rho^{-1}$  as follows
$$\rho^{-1}=\left\{
                  \begin{array}{ll}
                    \rho^{-1}_{-}, & \hbox{$x<x(t) $,} \\
                    0, & \hbox{$x=x(t)$,} \\
                    \rho^{-1}_{+}, & \hbox{$x>x(t)$.}
                  \end{array}
                \right.
\eqno{(3.12)}$$
 We  assert that (3.11) is a delta shock wave solution of  (1.1) with (1.2) in the sense of distributions if it satisfies the following generalized Rankine-Hugoniot relation
$$
\left\{
     \begin{array}{ll}
       \frac{dx(t)}{dt}=v_{\delta}(t), \\
       \frac{d}{dt}\Big(\frac{h(t)}{\sqrt{1-v_{\delta}^{2}(t)/c^{2}}}\Big)=v_{\delta}(t)\Big [\frac{n}{\sqrt{1-v^{2}/c^{2}}}\Big]-\Big[\frac{nv}{\sqrt{1-v^{2}/c^{2}}}\Big], \\
       \frac{d}{dt}\Big(\frac{w(t) v_{\delta}(t)}{1-v_{\delta}^{2}(t)/c^2}\Big)=v_{\delta}(t) \Big[\frac{\big(-\frac{1}{\rho c^{2}}+\rho \big) v}{1-v^2/c^2}\Big]
       -\Big[\frac{\big(-\frac{1}{\rho c^{2}}+\rho \big)v^2}{1-v^2/c^2}-\frac{1}{\rho}\Big ],\\
       \frac{d}{dt}\Big(\frac{w(t)}{1-v_{\delta}^{2}(t)/c^2}\Big)=v_{\delta}(t) \Big[\frac{\big(-\frac{1}{\rho c^{2}}+\rho \big)v^2/c^2}{1-v^2/c^2}+\rho\Big]
       -\Big[\frac{\big(-\frac{1}{\rho c^{2}}+\rho \big) v}{1-v^2/c^2}\Big],
     \end{array}
   \right.\eqno (3.13)
$$
  where $[q]= q_{+}-q_{-}$, etc.
 In fact,  if the relation (3.13) holds,  then for any test function $\phi\in C^{\infty}_{0}([0,+\infty)\times R),$  by Green's  formulation,  we have
 $$ \int_{0}^{+\infty}\int_{-\infty}^{+\infty}\Big(\Big(\frac{n}{\sqrt{1-v^{2}/c^{2}}}\Big)\phi_{t}+\Big(\frac{nv}{\sqrt{1-v^{2}/c^{2}}}\Big)\phi_{x}\Big)dxdt$$
 $$=\int_{0}^{+\infty}\int_{-\infty}^{x(t)}\frac{n_{-}}{\sqrt{1-v_{-}^{2}/c^{2}}}\phi_{t}+ \frac{n_{-}v_{-}}{\sqrt{1-v_{-}^{2}/c^{2}}} \phi_{x}dxdt$$
$$+\int_{0}^{+\infty}\int_{x(t)}^{+\infty}\frac{n_{+}}{\sqrt{1-v_{+}^{2}/c^{2}}}\phi_{t}+ \frac{n_{+}v_{+}}{\sqrt{1-v_{+}^{2}/c^{2}}} \phi_{x}dxdt$$ $$+\int_{0}^{+\infty}\frac{h(t)}{\sqrt{1-v_{\delta}^{2}(t)/c^{2}}}\phi_{t}(t,x(t))+ \frac{h(t)v_{\delta}(t)}{\sqrt{1-v_{\delta}^{2}(t)/c^{2}}} \phi_{x}(t, x(t))dt.\eqno{(3.14)}$$ \vskip 0.1in Without loss
of generality, we assume that $v_{\delta}(t):=\sigma_{\delta}$ is a constant and $\sigma_{\delta}>0.$  By exchanging  the ordering of integral and using
the change of variables,  the first term
on the right-hand side of (3.14) equals
$$\hspace{0mm}\int_{0}^{+\infty}\int_{-\infty}^{0}\frac{n_{-}}{\sqrt{1-v_{-}^{2}/c^{2}}}\phi_{t}dxdt
+ \int_{0}^{+\infty}\int_{0}^{x(t)}\frac{n_{-}}{\sqrt{1-v_{-}^{2}/c^{2}}}\phi_{t}dxdt$$
$$+\int_{0}^{+\infty}\int_{-\infty}^{x(t)}
\frac{n_{-}v_{-}}{\sqrt{1-v_{-}^{2}/c^{2}}} \phi_{x}dxdt
$$
$$\hspace{0mm}=\int_{0}^{+\infty}dx\int_{t(x)}^{+\infty}\frac{n_{-}}{\sqrt{1-v_{-}^{2}/c^{2}}}\phi_{t}dt
+\int_{0}^{+\infty}\frac{n_{-}v_{-}}{\sqrt{1-v_{-}^{2}/c^{2}}}\phi(t,x(t))dt$$$$
=-\int_{0}^{+\infty}\frac{n_{-}}{\sqrt{1-v_{-}^{2}/c^{2}}}\phi(t(x),x)dx
+\int_{0}^{+\infty}\frac{n_{-}v_{-}}{\sqrt{1-v_{-}^{2}/c^{2}}}\phi(t,x(t))dt$$
$$
=\int_{0}^{+\infty}\bigg(\frac{n_{-}v_{-}}{\sqrt{1-v_{-}^{2}/c^{2}}}-\frac{n_{-}v_{\delta}(t)}{\sqrt{1-v_{-}^{2}/c^{2}}}\bigg)\phi(t,x(t))dt.\eqno{(3.15)}$$
Similarly, the second term on the right-hand side of (3.14) equals
$$\hspace{0mm}\int_{0}^{+\infty}dx\int^{t(x)}_{0}\frac{n_{+}}{\sqrt{1-v_{+}^{2}/c^{2}}}\phi_{t}dt
-\int_{0}^{+\infty}\frac{n_{+}v_{+}}{\sqrt{1-v_{+}^{2}/c^{2}}} \phi(t, x(t))dt$$$$\hspace{0mm}
=\int_{0}^{+\infty}\frac{n_{+}}{\sqrt{1-v_{+}^{2}/c^{2}}}\phi(t(x),x)dx-\int_{0}^{+\infty}\frac{n_{+}v_{+}}{\sqrt{1-v_{+}^{2}/c^{2}}} \phi(t, x(t))dt$$
$$\hspace{0mm}
=\int_{0}^{+\infty}\bigg(\frac{n_{+}v_{\delta}(t)}{\sqrt{1-v_{+}^{2}/c^{2}}}-\frac{n_{+}v_{+}}{\sqrt{1-v_{+}^{2}/c^{2}}}\bigg ) \phi(t, x(t))dt.\eqno{(3.16)}$$

By using integrating by parts,  from (3.14)-(3.16), we obtain
$$ \int_{0}^{+\infty}\int_{-\infty}^{+\infty}\Big(\Big(\frac{n}{\sqrt{1-v^{2}/c^{2}}}\Big)\phi_{t}+\Big(\frac{nv}{\sqrt{1-v^{2}/c^{2}}}\Big)\phi_{x}\Big)dxdt
       $$
 $$=\int_{0}^{+\infty}\bigg(v_{\delta}(t)\Big [\frac{n}{\sqrt{1-v^{2}/c^{2}}}\Big]-\Big[\frac{nv}{\sqrt{1-v^{2}/c^{2}}}\Big] -
       \frac{d}{dt}\Big(\frac{h(t)}{\sqrt{1-v_{\delta}^{2}(t)/c^{2}}}\Big)\bigg)\phi(t, x(t))dt=0,\eqno{(3.17)}$$
which yields the first equality of (3.8).
Similarly, one can prove the second and third equalities of (3.8).  Thus, the assertion  is true.

\vskip 0.1in
\vskip 0.1in

\noindent{\small {\small\bf Remark 2.}  The generalized Rankine-Hugoniot relation (3.13)  reflects the exact relationship among the limit states on both sides of the delta shock wave and the  location, propagation speed, weights and the reassignment of $ v $ on the delta shock wave.

\indent
In addition,
  to guarantee uniqueness, we should propose the following entropy condition $$  \frac{v_{+}+\frac{1}{\rho_{+}}}{1+\frac{v_{+}}{\rho_{+} c^{2}}}
 \leq v_{\delta}(t)\leq  \frac{v_{-}-\frac{1}{\rho_{-}}}{1-\frac{v_{-}}{\rho_{-} c^{2}}}, \eqno{(3.18)}$$
which means that all characteristic  lines on both sides of  the delta shock wave are incoming.
A discontinuity  satisfying (3.13)  and (3.18) will be called a delta shock wave to system (1.1).

  At this moment,  the Riemann problem   is reduced to solve the
 ordinary differential equations (3.13) with  initial data
$$t=0:\,\,\,\,x(0)=0,\,\,\,v_{\delta}(0)=v_{0},\,\,\,h(0)=0,\,\,\,w(0)=0.\eqno{(3.19)}$$
where $v_{0}$ is an undetermined constant.

Integrating (3.13) from $0$ to $t$ with initial data (3.19), we have
$$\left\{\begin{array}{ll}
 \frac{h(t)}{\sqrt{1-v_{\delta}^{2}(t)/c^{2}}}=\Big [\frac{n}{\sqrt{1-v^{2}/c^{2}}}\Big]x(t)-\Big[\frac{nv}{\sqrt{1-v^{2}/c^{2}}}\Big]t, \\
       \frac{w(t) v_{\delta}(t)}{1-v_{\delta}^{2}(t)/c^2}= Fx(t)-Gt,\\
       \frac{w(t)}{1-v_{\delta}^{2}(t)/c^2}= Ex(t)-Ft,
     \end{array}
   \right.\eqno (3.20)
$$
   where
   $$E=\bigg[\frac{\big(-\frac{1}{\rho c^{2}}+\rho \big)v^2/c^2}{1-v^2/c^2}+\rho\bigg],F=\bigg[\frac{\big(-\frac{1}{\rho c^{2}}+\rho \big) v}{1-v^2/c^2}\bigg], G=\bigg[\frac{\big(-\frac{1}{\rho c^{2}}+\rho \big)v^2}{1-v^2/c^2}-\frac{1}{\rho}\bigg].$$

In what follows, under the  entropy condition (3.18), we can solve (3.20) to obtain
$$
\left\{
     \begin{array}{ll}
       x_{}(t)=\frac{F+\sqrt{F^{2}-EG}}{E}t, \\v_{\delta}(t)=\frac{F+\sqrt{F^{2}-EG}}{E},
       \\
       w(t) =\sqrt{F^{2}-EG}\bigg(1-\Big(\frac{F+\sqrt{F^{2}-EG}}{cE}\Big)^2\bigg)t,\\
       h(t)=
      \sqrt{1-\Big(\frac{F+\sqrt{F^{2}-EG}}{cE}\Big)^2}\,\bigg(\Big [\frac{n}{\sqrt{1-v^{2}/c^{2}}}\Big]\frac{F+\sqrt{F^{2}-EG}}{E}
     -\Big[\frac{nv}{\sqrt{1-v^{2}/c^{2}}}\Big]\bigg)t,
     \end{array}
   \right.\eqno (3.21)
$$
for $E\neq 0$, and

$$
\left\{
     \begin{array}{ll}
       x(t)=\frac{G}{2F}t, \\v_{\delta}(t)=\frac{G}{2F},
       \\
       w(t) =-F\bigg(1-\Big(\frac{G}{2Fc}\Big)^{2}\bigg )t,\\
       h(t)=\sqrt{1-\Big(\frac{G}{2Fc}\Big)^{2}}\,\bigg(\Big [\frac{n}{\sqrt{1-v^{2}/c^{2}}}\Big]\frac{G}{2F}-\Big[\frac{nv}{\sqrt{1-v^{2}/c^{2}}}\Big]\bigg)t,
     \end{array}
   \right.\eqno (3.22)
$$
for $E=0$.   The proof is similar to that in [10], so we omit it.

\vskip 0.1in
Thus, we have proved the following result.

\vskip 0.1in \noindent{\small {\small\bf Theorem 3.1.} On the physically relevant region,  under the condition $\frac{v_{+}+\frac{1}{\rho_{+}}}{1+\frac{v_{+}}{\rho_{+} c^{2}}}\leq\frac{v_{-}-\frac{1}{\rho_{-}}}{1-\frac{v_{-}}{\rho_{-} c^{2}}}$,
  Riemann problem  (1.1), (1.2) and (1.3)  admits  a unique  entropy solution in the sense of distributions of the form

$$(n, \rho,v)(t, x)=\left\{
                  \begin{array}{ll}
                    (n_{-}, \rho_{-}, v_{-}), & \hbox{$x<x(t) $,} \\
                    (h(t)\delta(x-x( t)),w(t)\delta(x-x( t)),   v_{\delta}(t)), & \hbox{$x=x(t)$,} \\
                    (n_{+}, \rho_{+}, v_{+}), & \hbox{$x>x(t)$,}
                  \end{array}
                \right.
\eqno{(3.23)}$$
where   $x(t)$, $v_{\delta}(t), $    $ h(t)$  and $w(t)$  are shown in (3.21) for $E\neq 0$ or
(3.22) for $E=0$.

\vskip 0.1in
At last, combining with the results in
Section 2, we  can conclude

\vskip 0.1in \noindent{\small {\small\bf Theorem 3.2.}  For
Riemann problem (1.1), (1.2) and (1.3),  on the physically relevant region,  there exists a unique entropy
solution, which consists of three contact discontinuities
when $\frac{v_{+}+\frac{1}{\rho_{+}}}{1+\frac{v_{+}}{\rho_{+} c^{2}}}>\frac{v_{-}-\frac{1}{\rho_{-}}}{1-\frac{v_{-}}{\rho_{-} c^{2}}}$
      and  a delta shock wave on which  both $\rho$ and $n$ contain  Dirac delta function   simultaneously when $\frac{v_{+}+\frac{1}{\rho_{+}}}{1+\frac{v_{+}}{\rho_{+} c^{2}}}\leq\frac{v_{-}-\frac{1}{\rho_{-}}}{1-\frac{v_{-}}{\rho_{-} c^{2}}}$.

\end{document}